# System Availability Optimization: Integrating Quantity Discounts and Delivery Lead Time Considerations

Zahra Sobhani*[1], Mahmoud Shahrokhi[2]

[1] Department of Industrial and Systems Engineering, University of Missouri, Columbia, MO, USA
[2] Department of Industrial Engineering, University of Kurdistan, Sanandaj, Iran

**Abstract**

The model allocates the system components orders to the suppliers to minimize the parts price and the system construction delay penalties and maximize the system availability during its use. It considers the quantity-based discount and variation of delivery lead time by ordering similar components. The model also reflects the prerequisite relationships between construction activities and calculates the delay penalty resulting from parts delivery lead time. This research presents a model for selecting suppliers of components of an industrial series-parallel multi-state system. A nonlinear binary mathematical program uses the Markov process results to select system components. It minimizes the total system construction phase costs, including the components' price, and the system construction delay penalty, and the system exploitation phase costs, including the system shutdown and working at half capacity. The model allocates the optimal orders for a typical industrial system's components, composing four elements. The proposed approach combines the nonlinear binary program and the Markov process results to optimize the system life cycle parameters, including the system construction cost and operational availability. Using the Markov chain results in binary nonlinear mathematical programming, this study attempts to strike the right balance between the construction phase's objectives and an industrial unit's operation phase.

---

[1] Corresponding author
Email address: zahra7sobhani@gmail.com

## 1. Introduction

When designing a system, there are always constraints related to the available manufacturing budget and time. On the other hand, system availability impacts the expected system costs and risks during the exploitation phase. This paper proposes an approach for supplier selection of the system components considering the purchasing price, delivery lead time, and parts' quality offered by different suppliers.

High-quality components improve the entire system's availability; however, they are more expensive. Lead times of suppliers may be different. Also, large orders benefit from the quantity-based discount; however, it leads to higher lead time because suppliers will need more time to prepare the additional components. The extended parts delivery lead times may cause delays in the system construction and result in a delayed penalty. This paper aims to minimize the product's life cycle cost, including the cost of building and operation, by considering the budget constraints and a minimum required reliability level. This paper presents a model to optimize the feedwater system of heat recovery steam generator (HRSG) boilers used in combined cycle power plants. Due to the non-negligible repair time, this paper uses the results of the continuous Markov chain approach in a nonlinear binary program. The following section presents a literature review on reliability design-based optimization studies. Then, after defining the problem, a simplified model of optimizing the feedwater system of HRSG boilers is developed. The final section interprets and discusses the results of the numerical example are.

## 2. Literature review

Many researchers applied availability and reliability analysis to optimize system configuration. Pascual and et al. present a model and a heuristic procedure to optimize the equipment's availability, procurement costs, and service levels subject to a budget constraint (Pascual et al., 2017). The model uses a single-echelon structure with repairable parts to determine the optimal number of repairable stocked spare parts. Results show that using an optimal procurement policy of spare parts produces better availability.

Continuous-time Markov Chain is a stochastic tool to describe and analyze systems behaviors under normal and failures conditions. It generates equations representing the portion of time the system stays in each state (Vasconcelos, 2017).

Billinton and Lee addressed an approach for reliability modeling utilizing the Markov approach for evaluating the design of a heat transport pump (Billinton and Lee, 1975). Using a continuous-time Markov chain model, Kumar et al. analyzed the chain-based scheme's availability and reliability (Kumar et al., 2009). Markov chain-based availability has been widely used in reliability-based design optimization (RBDO) to model system performance under uncertainty degradation and failure. This approach helps the analytical estimation of the steady-state or transition probability between system states, thereby calculating the proportion of time a system stays in operational, degraded, or failed states. Prior studies have illustrated the effectiveness of Markov models in showing component dependencies, repairs and replacements, and state-transition uncertainty, which make them particularly suitable for complex industrial systems (Sobhani & Shahrokhi, 2017a; Sobhani & Shahrokhi, 2017b; Sobhani, 2017c; Sobhani & Shahrokhi, 2019; Shahrokhi & Sobhani, 2018; Shahrokhi et al., 2021; Sobhani et al., 2025; Sobhani et al., 2023; Gharegozlu et al., 2018). These works highlight how availability-based performance metrics derived from Markov processes can be seamlessly integrated into optimization frameworks to support design decisions that balance reliability, cost, and system lifecycle considerations.

Che et al. introduced a Markov Chain-based availability model for studying the availability of virtual cluster systems. This model provides the possibility of explaining the lifecycle state and state transition of the virtual cluster node. They also obtained the availability level of each virtual cluster node in a complicated virtual cluster system or cloud data center (Che et al., 2011; Basiri et al., 2024). Musa Abbas et al. estimated the performance of a degraded multi-state system with three failure states using recursive Markov chain closed-form and an analytical solution considering a condition-based maintenance model (Musa Abbas et al., 2015). Huang et al. employed a Markov-chain-based model to assess the availability of wind turbines. They used a Poisson-process-based algorithm to estimate systems' state transition rates and solved a transition matrix to calculate expected availability (Huang et al., 2017). Siddiqui et al. optimized the availability of a multi-state series system configuration using the Markov model. They developed the optimal maintenance plan and validated the Monte Carlo Simulation results (Siddiqui et al., 2017).

Khojaste et al. applied Markov Decision Process models to optimize operational availability in energy systems, including electricity grids and wind-driven generation, demonstrating the use of stochastic state transitions to improve reliability and cost efficiency (Khojaste et al., 2025; Khojaste et al., 2024, see Khojaste et al., 2026 for a risk-averse version). Mooren addressed a new approach to making preliminary decisions for systems design. He utilized the Markov model and a network-theory-based approach to maximize the system availability of a naval vessel (Mooren, 2018).

Kumar and et al. used a semi-Markov process model to analyze the steady-state availability of mechanical systems that follow condition-based maintenance and evaluate optimal condition monitoring intervals (Kumar *et al.,* 2018). They used this methodology to analyze maintenance policies for a centrifugal pump. Tang et al. benefited from Markov Chain concepts in a service-based network model to design an availability evaluation approach for complex networks. Their method determines the minimum performance cost based on the quality of service metrics (Tang et al., 2020).

Afsharnia et al. used the Markov chain method to evaluate the availability of sugarcane harvester machines in the agro-industries (Afsharnia et al., 2020).

Sun et al. proposed a model to describe correlated fading processes in dual-frequency global navigation satellite system signals. They applied a Markov chain with stationary conditional probabilities based on the Poisson rate parameters by representing transitions in four fading states (Sun et al., 2021).

Traditional supplier selection methods often consider the purchase and shipping costs, ignoring the reliability and maintainability of purchased parts (Kanagaraj et al., 2014).

Jiao et al. proposed a new supplier house of quality approach consisting of customer requirement and supplier reliability designability. They also addressed the mathematic program model of supplier quality house to resolve the synthesis parts weighted scale (Jiao et al., 2008).

Kanagaraj et al. offered a reliability-based model to minimize procurement, maintenance, and downtime costs subjected to the practical constraints on product reliability and weight limitation and solved the model by simulated annealing algorithm (Kanagaraj and Jawahar, 2009). Tomaru et al. introduced new measures for supplier quality evaluation and a systematic approach to select suppliers. (Tomaru et al., 2011).

Kanagaraj et al. proposed a model for supplier selection problem using the reliability-based total cost of ownership approach regarding reliability and cost factors, including the product's initial price, quality, service, and maintenance-related costs. They also developed a meta-heuristic optimization algorithm by combining cuckoo search with a genetic algorithm to solve the problem (Kanagaraj et al., 2014). Hague et al. considered the selection of suppliers based on their ability to provide parts that improve the system's availability. They employed the availability of individual components as a criterion for the decision problem and weighted those criteria based on the component importance, and described the reliability and maintainability parameters by interval numbers (Hague et al., 2015). Yang also worked on the supplier selection problem. He employed operational cost and reliability factors to assess a supplier and then solved the model using a quantitative approach (Yang, 2016).

Rajabi Asadabadi presented a Markov chain to find a pattern for the changing preferences of customer needs. A combination of the analytic network process and quality function deployment links this pattern to product requirements and product requirements to supplier qualifications (Rajabi Asadabadi, 2017).

Azarian et al. explained eight essential practices, including reliability requirements and planning, training and development, reliability analysis, reliability testing, supply chain management, failure data tracking and analysis, verification and validation, and reliability improvement for supplier assessment problems (Azarian et al., 2020).

Some researchers studied the effects of the order quantities on the components parts and used them in supplier selection problems. Qin et al. considered all-units quantity-based discount in a freight allocation problem. They used a heuristic-based algorithm combining a filter-and-fan and a tabu search mechanisms (Qin et al., 2012). Also, Hamdan and Cheaitou worked on a multi-period green supplier selection and order allocation problem. They developed a single-product bi-objective integer linear programming model with deterministic demand considering all-unit quantity discounts was introduced (Hamdan and Cheaitou, 2017).

A literature review shows that researchers did not consider optimizing systems availability regarding components quantity-based discount and delivery lead time.

This paper presents a model to construct a specific series-parallel multi-state industrial system concerning availability, purchase price, and system construction delay cost. A binary mathematical program optimizes the system costs during its construction and exploitation phases simultaneously.

## 3. Problem statement

The studied problem concerns selecting components of some part of the boilers' feedwater system used in power plants. It is a multi-state parallel series system, including subsystems shown in figures 1 and 2.

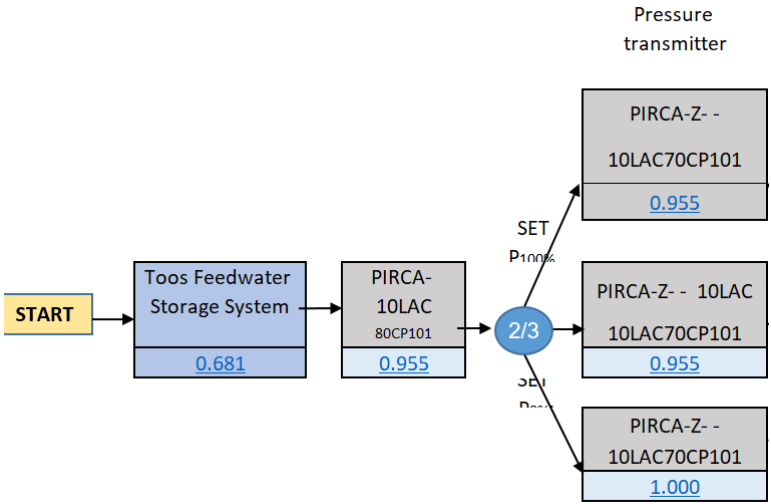

Figure 1. The reliability block diagram of the studied series-parallel system

Generally, The suppliers of two serial components are identical. Combining these components to form part A, in figure 2, simplifies the model.

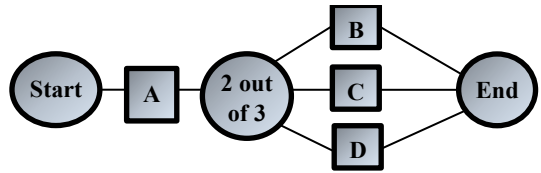

**Figure 2. Simplified studied series-parallel system**

Several alternative suppliers offer different prices, qualities (failure and repair rates), and delivery lead times. The availability of the entire system depends on the reliability and repair rate of its subsystems. Components' reliability and repair rate, price, and delivery lead time are criteria used to allocate orders to the suppliers.
Despite quantity-based discounts, large orders increase the delivery lead time when ordering several parts from a single supplier.

Each of B, C, and D can provide 50% of the nominal capacity. Therefore, the system states are: (1) full capacity (FC), when A operates and at least two out of B, C, and D are operational, (2) half-capacity (HC), when component A operates, and one out of B, C and D are working, and (3) shutdown (SD), when A or all B, C, and D are not functional. Occurring HC and SD cause lost production and repair costs.

On the one hand, the higher reliable components increase the entire system's availability; on the other hand, they are generally more expensive and have a longer delivery lead time. Besides, there is a daily penalty for system delays in system completion time. The model's objective is to minimize the total cost of system assembly and exploitation.

The Markov chain model represents a class of stochastic processes in which the probability of future events depends only on the present state (Johnson, 2020). When the process evolves, the system can remain in the same state or change to a different state. Transition probabilities describe the likelihood of movements between states. As the following equation shows, the transition matrix (Markov matrix) for a system having M states indicates all the transition probabilities ($P_{ij}$), suggesting the likelihood of changing the system state from i to j.

$$P = \begin{bmatrix} p_{11} & p_{12} & \cdots & p_{1M} \\ p_{21} & p_{22} & \cdots & p_{2M} \\ \cdots & \cdots & \cdots & \cdots \\ p_{M1} & p_{M2} & \cdots & p_{MM} \end{bmatrix} \quad (1)$$

All the matrix elements are nonnegative, and the sum of the numbers in each row equals 1. A state-space diagram shows system states, transition paths, and related transfer rates. Figure 4 shows the state space diagram for the feedwater system used in HRSG boilers. The nodes in this diagram represent the system's state, and the branches represent the transition from one state to another. The values on the lines indicate the change rate of moving from one state to another. It may signify λ or μ, meaning the system failure and repair rates, respectively (Mohammadi et al., 2021).

The following section presents the model assumptions, parameters, decision variables, objective function, and constraints.

### 3.1. Assumptions

This study uses a continuous-time Markov chain to obtain the probability of each state, where the time that system spends in each state before proceeding to the next state are independent exponential random variables (Ross, 2019). The other assumptions of the model are as follows:

1. The construction project has a deadline, and the delay will impose a fixed daily cost to the project.
2. There are several alternative suppliers For each component.
3. The time and cost of components assembly processes are independent of the quality of the components.
4. Suppliers offer components with known quality (reliability and repair time).
5. Supplier delivery lead time and price depend on the order's quantity.
6. B, C, and D functions in the system are identical, and each provides 50% of the nominal system capacity.
7. The system always will be in one of the three states: (1) full capacity (FC), (2) half capacity HC) and (3) shutdown (SD), according to failure or working of its components
8. There are daily loss production penalties for HC and SD states.
9. The components' mean time to repair is not negligible.
10. The system has an infinite lifetime.
11. Failure and repair rates of the components are constant during the planning horizons.
12. The purchase price includes the cost of ordering, shipment, and assembly of the components.
13. The construction of B, C, and D should begin after completing the assembly process of A.
14. The economic minimum attractive return rate (MARR) is known.

### 3.2. Indices

$i = 1,2,3,4$     The components' index (i=1, 2, 3, and 4 for A, B, C, and D, respectively)
$j = 1,2,...,J$     The suppliers' index
$k = 1,2,...,K$     Index of states of the system

### 3.3. Parameters

$C_{1j}$:    The purchase price of component A from supplier j
$C_{ij}$:    The unit purchase price of B, C, or D from supplier j, when ordered separately
$C'_{ij}$:    Unit purchase price when two of components B, C, or D ordered to supplier j
$C''_{ij}$:    Unit purchase price when B, C, and D ordered to supplier j
$R_{ij}$:    Reliability of component i from supplier j

$R$:    The requested minimum system's reliability
$T$:    System construction deadline
$B$:    The maximum available budget for purchasing components
$C_2$:    The daily cost of SD
$C_3$:    The daily cost of working the system at half capacity

| | |
|---|---|
| $C^D$: | Daily delay penalty for the system construction |
| $L_{Aj}$: | The delivery lead time of component A from supplier j |
| $L_{Bj}$: | The delivery lead time of components B, C, or D from supplier j, when ordered separately |
| $L'_{Bj}$: | The delivery lead time, when two components of B, C, or D ordered to supplier j |
| $L''_{Bj}$: | The delivery lead time when all B, C, and D, ordered to supplier j |
| $F_{ik}$: | Time of the process k on component i |
| $K_1$: | Set of the required processes for construction component 1 |
| $K_2$: | Set of the required processes for construction components 2, 3 and 4 |
| r: | The minimum acceptable rate of return (MARR) |
| $\alpha$: | The hourly cost of SD. |
| $\beta$: | The hourly cost of HC. |
| $\mu_{ij}$: | Repair rate of component i from supplier j expressed in terms of the average number of performed repairs actions per hour |
| $\lambda_{ij}$: | The failure rate of component i from supplier j (expressed in terms of the average number of failures per hour) |
| $A_1$: | Total purchase cost |
| $A_2$: | The expected total cost of the system shutting down (in the entire life cycle of the system) |
| $A_3$: | The total cost of the working system at half capacity (in the entire life of the system) |
| $A_4$: | Total system construction delay penalty |

### 3.4. Decision Variables

| | |
|---|---|
| $y_{ij}$: | A binary variable; equal to 1 if component i is purchased from supplier j and zero otherwise |
| $R_i$: | Reliability of component i |
| $R_e$: | Reliability of the entire system |
| $P_0$: | The proportion of the time that the system is shutting down |
| $P_{50}$: | The proportion of the time that the system is working at half capacity |
| $S_k$: | The proportion of the time that the system is in state k |
| $\mu_i$: | Repair rate of component i (expressed in terms of the average number of repairs actions performed per hour) |
| $\lambda_i$: | The failure rate of component i (expressed in terms of the average number of failures per hour) |
| $T^c$: | Completion time of the system construction project |
| $T_c^i$: | Completion time of construction component i |
| L: | Required time for delivering all components B, C, and D |

## 4. Mathematical model

This study develops an approach for determining the optimal configuration for a typical industrial system. The same method is applicable for other cases. The binary nonlinear programming (BLP) model is as follows. The model's objective is to minimize the system's expected total cost, including purchase cost, project construction delay penalty, and total cost of being system in shutdown and half capacity states during its exploitation.

$$Min\ Z = A_1 + A_2 + A_3 + A_4 \qquad (1)$$

Constraint (2) uses binary variables to calculate purchase costs. It equals the sum of the purchase price of all components regarding the quantity-based discounts.

$$A_1 = \sum_{j=1}^{3} 3 y_{2j} y_{3j} y_{4j} C_{ij}''$$

$$+ \sum_{j=1}^{3} \begin{pmatrix} y_{2j} y_{3j} (1-y_{4j}) \\ + y_{2j} y_{4j} (1-y_{3j}) \\ + y_{3j} y_{4j} (1-y_{2j}) \end{pmatrix} \times C_{ij}'$$

$$+ \sum_{j=1}^{3} \begin{pmatrix} y_{2j}(1-y_{3j})(1-y_{4j}) \\ + y_{3j}(1-y_{2j})(1-y_{4j}) \\ + y_{4j}(1-y_{2j})(1-y_{3j}) \end{pmatrix} \times C_{ij} \qquad (2)$$

$$+ \sum_{j=1}^{3} y_{1j} C_{1j}$$

Constraint (3) represents the total cost of system shut down during its infinite exploitation period by considering the capital rate of return.

$$A_2 = 8760 * \alpha P_0 * \frac{1}{r} \qquad (3)$$

Equation (4) calculates the system's proportional time in shutdown states when part A or all three components B, C, and D fail.

$$P_0 = S_8 + S_9 + S_{10} + S_{11} \\ + S_{12} + S_{13} + S_{14} + S_{15} \qquad (4)$$

Constraint (5) represents the total cost of the working system at half capacity (when the system works at 50% of its power), considering the capital rate of return during its infinite exploitation period.

$$A_3 = 8760 * \beta P_{50} * \frac{1}{r} \qquad (5)$$

Constraint (6) shows the delay penalty of the system construction project. If the plan completes before the deadline ($T^c \leq T$), the delay penalty equals zero.

$$A_4 = \max(0, T^c - T) C^D \qquad (6)$$

Equation (7) calculates the completion date of component A installation. It is the sum of its lead time and needed installation time, which is a constant amount.

$$T_c^1 = \sum_{j=1}^{3} y_{1j} L_{Aj} + \sum_{k \in k_1} F_{1k} \qquad (7)$$

Formula (8) calculates the completion date of installation of B, C, and D. If the component lead time (L) is bigger than the completion date of component A, then its completion date is the sum of its lead time and its needed installation time. In case of component A is installed after its lead time, $T_c^i$ equals completion date of component A plus installation time of part i.

$$T_c^i = \max(L, T_c^1) + \sum_{k \in k_2} F_{ik}$$

$$i = 2, 3 \text{ and } 4 \tag{8}$$

According to their suppliers, constraint (9) gives the delivery lead time for components B, C, and D.

$$L = \sum_{j=1}^{3} 3 y_{2j} y_{3j} y_{4j} L_{Bj}''$$
$$+ \sum_{j=1}^{3} \begin{pmatrix} y_{2j} y_{3j} (1 - y_{4j}) \\ + y_{2j} y_{4j} (1 - y_{3j}) \\ + y_{3j} y_{4j} (1 - y_{2j}) \end{pmatrix}$$
$$\times Max(L_{Bj}', L_{Bj})$$
$$+ \sum_{j=1}^{3} \begin{pmatrix} y_{2j} (1 - y_{3j})(1 - y_{4j}) \\ + y_{3j} (1 - y_{2j})(1 - y_{4j}) \\ + y_{4j} (1 - y_{2j})(1 - y_{3j}) \end{pmatrix}$$
$$\times L_{Bj} \tag{9}$$

Equation (10) represents the completion date of project construction.

$$T_c = \max(T_c^i) \tag{10}$$

Equation (11) calculates the system's proportional time in half capacity states (states 5, 6, and 7).

$$P_{50} = S_5 + S_6 + S_7 \tag{11}$$

Constraint (12) indicates the maximum available budget for buying components.

$$A_1 \leq B \tag{12}$$

Equations (13), (14), and (15) represent the reliability, repair rate, and failure rate of each element, respectively, according to the selected supplier.

$$R_i = \sum_{j=1}^{3} Y_{ij} R_{ij} \qquad \forall i \tag{13}$$

$$\mu_i = \sum_{j=1}^{3} Y_{ij} \mu_{ij} \qquad \forall i \tag{14}$$

$$\lambda_i = \sum_{j=1}^{3} Y_{ij} \lambda_{ij} \qquad \forall i \tag{15}$$

Constraint (16) ensures selecting only one supplier for each component.

$$y_{ij} \in \{0.1\} \qquad (16)$$

Equation 17 defines variables $y_{ij}$ as binary variables.

$$\sum_{j=1}^{3} y_{ij} = 1 \qquad \forall i \qquad (17)$$

Equations (18) to (32) are equilibrium equations related to the continuous Markov process model, which are written according to figure (4) and ensure that input streams are equal to the output streams.

$$(\lambda_1 + \lambda_2 + \lambda_3 + \lambda_4)S_1 = \mu_2 S_2 + \mu_4 S_3 + \mu_3 S_4 + \mu_1 S_{13} \qquad (18)$$

$$(\lambda_1 + \mu_2 + \lambda_3 + \lambda_4)S_2 = \lambda_2 S_1 + \mu_4 S_5 + \mu_3 S_6 + \mu_1 S_9 \qquad (19)$$

$$(\lambda_1 + \mu_4 + \lambda_3 + \lambda_2)S_3 = \lambda_4 S_1 + \mu_2 S_5 + \mu_1 S_{12} + \mu_3 S_7 \qquad (20)$$

$$(\lambda_1 + \mu_3 + \lambda_2 + \lambda_4)S_4 = \lambda_3 S_1 + \mu_4 S_7 + \mu_2 S_6 + \mu_1 S_{14} \qquad (21)$$

$$(\lambda_1 + \mu_2 + \lambda_3 + \mu_4)S_5 = \lambda_4 S_2 + \lambda_2 S_3 + \mu_3 S_{11} + \mu_1 S_8 \qquad (22)$$

$$(\lambda_1 + \mu_2 + \lambda_4 + \mu_3)S_6 = \lambda_3 S_2 + \lambda_2 S_4 + \mu_4 S_{11} + \mu_1 S_{10} \qquad (23)$$

$$(\lambda_1 + \mu_3 + \lambda_2 + \mu_4)S_7 = \lambda_4 S_4 + \lambda_3 S_3 + \mu_2 S_{11} + \mu_1 S_{15} \qquad (24)$$

$$\mu_1 S_8 = \lambda_1 S_5 \qquad (25)$$

$$\mu_1 S_9 = \lambda_1 S_2 \qquad (26)$$

$$\mu_1 S_{10} = \lambda_1 S_6 \qquad (27)$$

$$(\mu_2 + \mu_3 + \mu_4)S_{11} = \lambda_4 S_6 + \lambda_2 S_7 + \lambda_3 S_5 \qquad (28)$$

$$\mu_1 S_{12} = \lambda_1 S_3 \qquad (29)$$

$$\mu_1 S_{13} = \lambda_1 S_1 \qquad (30)$$

$$\mu_1 S_{14} = \lambda_1 S_4 \qquad (31)$$

$$\mu_1 S_{15} = \lambda_1 S_7 \qquad (32)$$

Constraint (33) ensures that the sum of all probabilities of being system in all states is equal to 1.

$$\sum_{i=1}^{15} S_i = 1 \qquad (33)$$

Constraint (34) specifies the reliability of the entire system.

$$R_e = 1 - P_0 \qquad (34)$$

Constraint (35) ensures that the reliability of the entire system is more than a predefined level.

$$R_e \geq R \qquad (35)$$

## 5. Numerical example

The presented Markov diagram in figure (4) shows the set of all possible states and their relationships. In this figure, oval, trapezoid, and rectangular forms represent FC, HC, and SD states.

Tables (I) shows values of the basic parameters of the numerical example.

Table I. The values of basic parameters of the numerical example

| Parameter | Value | Parameter | Value |
|---|---|---|---|
| $C2$ | 80000 | $\mu_{i2}$ | 0.07 |
| $C3$ | 30000 | $\mu_{i3}$ | 0.1 |
| $C^D$ | 300 | $\lambda_{i1}$ | 0.05 |
| $B$ | 1100 | $\lambda_{i2}$ | 0.03 |
| $T$ | 68 | $\lambda_{i3}$ | 0.01 |
| $\alpha$ | 0.2 | $R_{i1}$ | 0.9 |
| $\beta$ | 0.1 | $R_{i2}$ | 0.95 |
| $r$ | 0.1 | $R_{i3}$ | 0.99 |
| $\mu_{i1}$ | 0.05 | $R$ | 0.8 |
| i=1,2,3,4 | | | |

Table (II) illustrates the components price offered by different suppliers, according to their order quantity.

Table II. The price of components

| Component price | | Supplier 1 | Supplier 2 | Supplier 2 |
|---|---|---|---|---|
| $C_{1j}$ | | 200 | 220 | 240 |
| i=2,3,4 | $C_{ij}$ | 300 | 340 | 380 |
| | $C'_{ij}$ | 250 | 280 | 320 |
| | $C''_{ij}$ | 200 | 240 | 280 |

Table (III) illustrates the required time for the manufacturing processes.

Table III. Duration of manufacturing processes

| | $F_{i1}$ | $F_{i2}$ | $F_{i3}$ | $F_{i4}$ | $F_{i5}$ |
|---|---|---|---|---|---|
| $F_{1k}$ | 3 | 5 | 7 | 4 | 2 |
| $F_{ik}$ $i=2,3,4$ | 6 | 13 | 16 | 5 | - |

Table (IV) presents the delivery lead time of the components, according to their order quantity.

Table IV. The components' delivery lead time

|  | Supplier 1 | Supplier 2 | Supplier 2 |
|---|---|---|---|
| $L_{Aj}$ | 5 | 7 | 17 |
| $L_{ij}$ | 6 | 19 | 31 |
| $L'_{Bj}$ | 8 | 24 | 37 |
| $L''_{Bj}$ | 12 | 30 | 42 |

The second column of table (V) shows the optimal solution, including variables and their value. In the resolution, the reliability of components A is more than other components, and B and D have the same quality. The system's overall reliability is about 0.85, and the system shutdowns 15.3% of the time. It also works at half capacity 29.7% of the time. The project delay is 24 days, and the delay penalty is a decisive parameter and more than half of the total cost. The remaining expense is due to the failure of components and their purchase price.

## 6. Sensitivity analysis

This section explains the sensitivity of the optimal solution to the variation of the model parameters.

### 6.1. The effects of quantity-based discount

The third column of table V shows the result of resolving the model without quantity-based discount. The model has no feasible solution in this case, and under the lack of price reduction, the requested minimum system reliability (0.8) requires at least 1140 dollars to buy high-quality components. Even by considering 1200 dollars budget, the price of the high-quality components increases, and the model selects lower reliability elements, which increases the system failures cost (i.e., half capacity and shutdown costs) due to the reduction of entire system availability.

On the other hand, it reduces the project's completion time and consequently delays penalty by selecting components with shorter delivery lead times. In general, removing the quantity-based discounts cause an augmentation in the purchase and operating costs and a reduction in system project delay penalty.

### 6.2. The effects of delivery lead time

Column 4 of table (V) illustrates the results of resolving the model by a fixed delivery lead time, independent of the order sizes. The model selects all B, C, and D with similar high reliability (0.95) without a radical cost augmentation. Therefore, the overall availability improves, which reduces the HC and SD costs. Despite ordering components to suppliers whose delivery lead times are longer, the project completion time decreases. Total cost also decreases as a result of failure cost and delay penalty reduction.

## 6.3. The effects of available budget and the project delay penalty

As shown in column 5 of table (V), the optimal solution has not changed by increasing the budget to 1200 dollars. In fact, by allocating more funds to purchase components, the model does not select higher quality components to impede high delivery lead time. It means the reduction of HC and SD costs is less significant than the delay penalty increment.

According to Columns 6 and 7, a reduction in project delay daily penalty to 100 $/day does not affect the components' ordering plan and, consequently, the entire system's reliability and availability. Therefore, the cost of purchasing, working at half capacity, and the system shutdown do not change.

However, the components' quality has increased after reducing the daily project penalty to 50 $/day. In this case, the system's availability increases while the system failures cost reduces. Compared to the optimal solution, the project's completion time increases significantly (44 days), leading to a 1.5 times increment in the project delay time. Therefore, in this example, the slight reduction in the delay penalty rate does not significantly affect the components' order plan. After incrementing the delay penalty to 500 $/day, higher quality components are not economical anymore. So purchase and operation costs remain constant and only delay penalty increases.

Table V. The optimal solution including variables and their value

| Variable (1) | Optimal Solution (2) | Sensitivity Analyses | | | | | |
|---|---|---|---|---|---|---|---|
| | | Removing discount (3) | fixed delivery lead time (4) | B=1200 (5) | B=1200 & $C^D$=100 (6) | B=1200 and $C^D$=50 (7) | B=1200 & $C^D$=500 (8) |
| $R_1$ | 0.99 | 0.99 | 0.99 | 0.99 | 0.99 | 0.99 | 0.99 |
| $R_2$ | 0.9 | 0.9 | 0.95 | 0.9 | 0.95 | 0.95 | 0.9 |
| $R_3$ | 0.95 | 0.9 | 0.95 | 0.95 | 0.9 | 0.95 | 0.95 |
| $R_4$ | 0.9 | 0.9 | 0.95 | 0.9 | 0.9 | 0.99 | 0.9 |
| $Z$ | 13572 | 13284.1 | 7453.5 | 13572 | 8772.1 | 7247.6 | 18372 |
| $R_e$ | 0.847 | 0.805 | 0.887 | 0.847 | 0.847 | 0.9 | 0.847 |
| $T^c$ | 92 | 87 | 78 | 92 | 92 | 136 | 92 |
| $A_1$ | 1080 | 1140 | 960 | 1080 | 1080 | 1180 | 1080 |
| $A_2$ | 2686.13 | 3423.5 | 1984.6 | 2686.13 | 2686.13 | 1711.3 | 2686.13 |
| $A_3$ | 2605.95 | 3020.7 | 1508.8 | 2605.95 | 2605.95 | 956.35 | 2605.95 |
| $A_4$ | 7200 | 5700 | 3000 | 7200 | 2400 | 3400 | 12000 |
| $P_0$ | 0.153 | 0.195 | 0.113 | 0.153 | 0.153 | 0.098 | 0.153 |
| $P_{50}$ | 0.297 | 0.345 | 0.172 | 0.297 | 0.297 | 0.109 | 0.297 |
| $P_{100}$ | 0.549 | 0.46 | 0.714 | 0.549 | 0.549 | 0.793 | 0.549 |

# 7. Discussions

This paper considers a reliability allocation problem regarding supply chain parameters, including price, reliability, and lead time of the components. The proposed model is applicable in the three-dimensional concurrent engineering context, simultaneously considering component supply, system construction, and exploitation phases. The model allocates the orders of components of a typical multi-state serial-parallel system to suppliers. Developing a similar model for other systems requires redrawing the reliability block diagram and rewriting the Markov process equations. The proposed approach is applicable in various industries to improve the systems' availability and reduce the products' final price.

The model considers the following suppliers' financial, operational, and essential quality parameters during the system's design, construction, and exploitation:

- Among economic parameters, it includes the components price and quantity-based discount.
- Construction parameters are delivery lead times and components assembly time, affecting the system construction date and the project delay penalty.
- The operation parameters are components failure and repair rates, which affect HC and SD duration.

The model also supposed some assumptions:

- They knowen fixed components failure and repair rates during the production system's life cycle. In other words, the time between failure and repair times are random variables with exponential distribution. Many types of research used these assumptions, which are essential for developing the Markov model. Considering constant failure rate is a commonly used assumption for the mature phase of the system's exploitation.
- Considering components' repair rate as a Normal random variable improves the model. In this case, the application of the Markov process for this problem will require simplification assumptions.
- Considering infinite lifetime for the system is acceptable for the typical feedwater system of powerplants boilers. Viewing a known limited system lifetime requires adapting the engineering economic factor in the objective function. It does not change the generality of the model. However, some systems have a finite lifetime. So it is helpful to evaluate and better the model considering limited lifetime.
- HC and SD cost functions are independent of the duration of these states. Defining HC and SD cost functions regarding both HC and SD times and duration improves the model.

- Reliability, price, and lead time factors are constant in this paper. Future researches can enhance the model by considering uncertainty in components' reliability, cost, and lead times.

- Studying the effect of repair rates increment caused by the similarity of elements (i.e., batch orders) has not been investigated. It is also a suggestion for future researches.

## 8. Conclusions

This paper develops a nonlinear binary programming model to optimize a multi-state parallel-series system's construction and exploitation cost by determining the components' ordering plan. The model minimizes the system's total cost, including system components purchase price, system construction delay penalty, reduced system capacity, and shutdown cost. The model considers the effects of components' order size on their buying price and delivery lead time. In addition, the paper studied the sensitivity of the following parameters on the optimal solution:

1- quantity-based discount
2- variable delivery lead time
3- available budget
4- the daily penalty rate

Numerical example sensitivity analyses showed that:

1- Reducing the delay penalty rate to 50 $/day has the most significant effect on order planning and, as a result, all variables of the model.
2- Removing quantity-based discounts made the model choose elements with lower quality. At the same time, fixed delivery lead time and the delay penalty reduction to 50 $/day caused the selection of high-quality components.
3- Changing the available budget and the daily penalty rate to 100 and 500 did not affect the results.

The delay in the exampled system construction project significantly impacts its components orders. The components similarity discount reduces the required initial budget and increases the availability of the entire system. Results also show effects of variation in components delivery lead times, components similarity, components purchase budget on the total system cost.

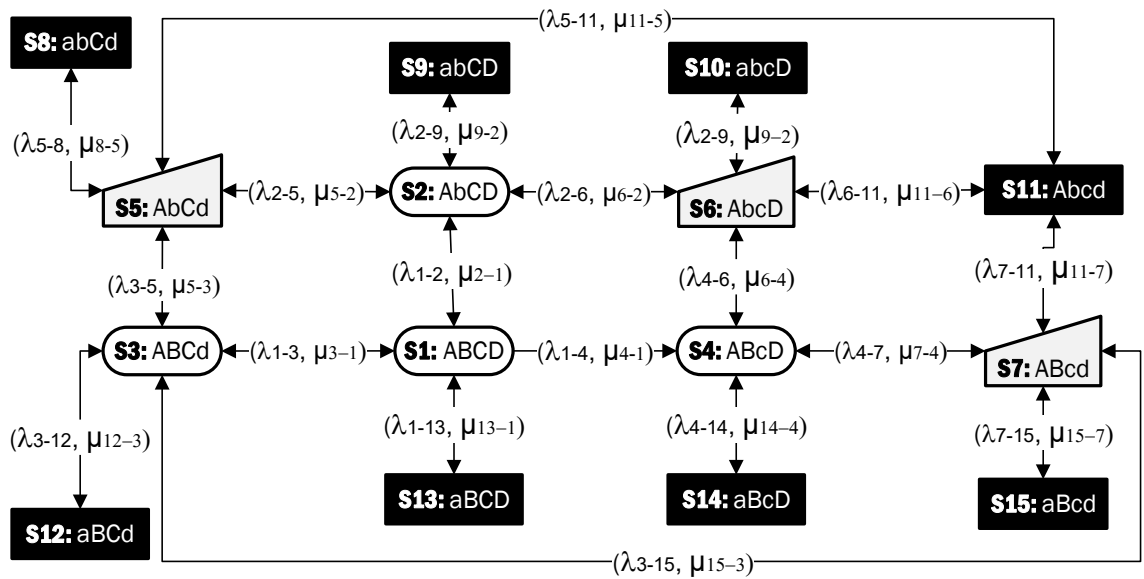

**Figure 4. The system Markov chain model**